\documentclass{birkjour}
\usepackage{mathtools}
 \newtheorem{thm}{Theorem}[section]
 
 \newtheorem{lem}[thm]{Lemma}
 
 \theoremstyle{definition}
 \newtheorem{defn}[thm]{Definition}
 \theoremstyle{remark}
 
 \newtheorem{ex}{Example}[section]
 \numberwithin{equation}{section}
\begin{document}

\title[Para-Sasakian metric as an almost conformal $\eta$-Ricci soliton]
{Geometry of para-Sasakian metric as an almost conformal $\eta$-Ricci soliton}

\author[S. Sarkar]{Sumanjit Sarkar}
\address{Department of Mathematics\\
Jadavpur University\\
Kolkata-700032, India.}
\email{imsumanjit@gmail.com}

\author[S. Dey]{Santu Dey}
\address{Department of Mathematics\\
Bidhan Chandra College\\
Asansol, Burdwan, West Bengal-713304, India.}
\email{santu.mathju@gmail.com}

\author[A. Bhattacharyya]{Arindam Bhattacharyya}
\address{Department of Mathematics\\
Jadavpur University\\
Kolkata-700032, India}
\email{bhattachar1968@yahoo.co.in}


\begin{abstract}
In this paper, we initiate the study of conformal $\eta$-Ricci soliton and almost conformal $\eta$-Ricci soliton within the framework of para-Sasakian manifold. We prove that if para-Sasakian mteric admits conformal $\eta$-Ricci soliton, then the manifold is $\eta$-Einstein and either the soliton vector field $V$ is Killing or it leaves $\phi$ invariant. Here, we have shown the characteristics of the soliton vector field $V$ and scalar curvature when the manifold admitting conformal $\eta$-Ricci soliton and vector field is pointwise collinear
with the characteristic vector field $\xi$. Next, we show that a para-Sasakian metric endowed an almost conformal $\eta$-Ricci soliton is $\eta$-Einstein metric if the soliton vector field $V$ is an infnitesimal contact transformation. We have also displayed that the manifold is Einstein if it represents a gradient almost conformal $\eta$-Ricci soliton. We have developed an example to display the alive of conformal $\eta$-Ricci soliton on 3-dimensional para-Sasakian manifold.
\end{abstract}

\maketitle

\section{\textbf{Introduction}}

In modern mathematics, the methods of paracontact geometry play an important role. The notion of almost paracontact manifold was first introduced by Sato \cite{Sato}. After that he and Matsumoto \cite{Matsumoto} defined and studied a para-Sasakian manifold as special case of an almost paracontact manifold. Adati et al. \cite{Adati} deduced some fundamental properties of para-Sasakian manifold. Later Kaneyuki and Williams \cite{Kaneyuki} associated pseudo-Riemannian metric with an almost paracontact manifold after Takahashi \cite{Takahashi} intoduced pseudo- Riemannian metric in contact manifold, in particular, in Sasakian manifold. Zamkovoy in \cite{Zamkovoy} proved that any almost paracontact structure admits a pseudo-Riemannian metric with signature $(n+1,n)$. Para-Sasakian manifold (in short p-Sasakian manifold) was studied by many authors, namely: Calvaruso \cite{Calvaruso-1}, Cappelletti et al. \cite{Cappelletti}, Tripathi et al. \cite{Tripathi} and many others.

A pseudo-Riemannian manifold $(M,g)$ admits a Ricci soliton which is a generalization of Einstein metric (i.e, $S=ag$ for some constant $a$) if there exists a smooth non-zero vector field $V$ and a constant $\lambda$ such that,
\begin{equation*}
  \frac{1}{2}\mathcal{L}_Vg+S+\lambda g=0,
\end{equation*}
where $\mathcal{L}_V$ denotes Lie derivative along the direction $V$ and $S$ denotes the Ricci curvature tensor of the manifold. The vector field $V$ is called potential vector field and $\lambda$ is called soliton constant.\par

The Ricci soliton is a self-similar solution of the Hamilton's Ricci flow \cite{Ham} which is defined by the equation $\frac{\partial g(t)}{\partial t}=-2S(g(t))$ with initial condition $g(0)=g$, where $g(t)$ is a one-parameter
family of metrices on $M$. The potential vector field $V$ and soliton constant $\lambda$ play  vital roles while determining the nature of the soliton. A soliton is said to be shrinking, steady or expanding according as $\lambda<0$,
$\lambda=0$ or $\lambda>0$. Now if $V$ is zero or Killing then the Ricci soliton reduces to Einstein manifold and the soliton is called trivial soliton.\par

If the potential vector field $V$ is the gradient of a smooth function $f$, denoted by $Df$ then the soliton equation reduces to,
\begin{equation*}
  Hessf+S+\lambda g=0,
\end{equation*}
where $Hessf$ is Hessian of the smooth function $f$. Perelman \cite{Perelman} proved that a Ricci soliton on a compact manifold is a gradient Ricci soliton.\par

In 2005, Fischer \cite{Fischer} has introduced conformal Ricci flow which is a mere generalisation of the classical Ricci flow equation that modifies the unit volume constraint to a scalar curvature constraint. The conformal Ricci flow equation was given by,
\begin{eqnarray*}
  \frac{\partial g}{\partial t}+2(S+\frac{g}{n}) &=& -pg, \\
  r(g) &=& -1,
\end{eqnarray*}
where $r(g)$ is the scalar curvature of the manifold, $p$ is scalar non-dynamical field and $n$ is the dimension of the manifold. Corresponding to the conformal Ricci flow equation in 2015, Basu and Bhattacharyya \cite{Basu} introduced the notion of conformal Ricci soliton equation as a generalization of Ricci soliton equation given by,
\begin{equation*}
\mathcal{L}_Vg+2S+[2\lambda-(p+\frac{2}{n})]g=0.
\end{equation*}

In 2009, Cho and Kimura \cite{Kimura} introduced the concept of $\eta$-Ricci soliton which is another generalization of classical Ricci soliton  and is given by,
\begin{equation*}
  \mathcal{L}_\xi g+2S+2\lambda g+2\mu\eta\otimes\eta=0,
\end{equation*}
where $\mu$ is a real constant, $\eta$ is a 1-form defined as $\eta(X)=g(X,\xi)$ for any $X\in\chi(M)$. Clearly it can be noted that if $\mu=0$ then the $\eta$-Ricci soliton reduces to Ricci soliton.\par

Recently Siddiqi \cite{Siddiqi} established the notion of conformal $\eta$-Ricci soliton which generalizes both conformal Ricci soliton and $\eta$-Ricci soliton. The equation for conformal $\eta$-Ricci soliton is given by,
\begin{equation}\label{A1}
\mathcal{L}_Vg+2S+[2\lambda-(p+\frac{2}{n})]g+2\mu\eta\otimes\eta=0.
\end{equation}
In the foregoing equation if we consider the soliton vector field as a gradient of a smooth function $f$ and $\mu$ as a smooth function then the soliton equation changes to
\begin{equation}\label{A2}
  Hess f+S+[\lambda-(\frac{p}{2}+\frac{1}{(2n+1)})]g+\mu\eta\otimes\eta=0,
\end{equation}
and the soliton is called gradient almost conformal $\eta$-Ricci soliton.\\
As follows in the literature, Ricci soliton on paracontact geometry studied by many authors (\cite{Bejan}, \cite{Calvaruso}, \cite{Pra}). In particular, Calvaruso and Perrone \cite{Calvaruso} explicitly studied Ricci soliton on 3-dimensional almost paracontact manifolds. In 2019, Patra \cite{Patra} studied Ricci soliton on paracontact metric manifolds and proved that if a para-Sasakian manifold satisfy Ricci soliton equation then the manifold is either Einstein or $\eta$-Einstein. The case of $\eta$-Ricci soliton in para-Sasakian manifold was treated by Naik and Venkatesha in \cite{Naik} and showed that if the metric of a para-Sasakian manifold represents a $\eta$-Ricci soliton then the manifold is either Einstein or D-homothetically fixed $\eta$-Einstein manifold. Very recently conformal $\eta$-Ricci soliton and its generalizations have been studied by \cite{roy, sarkar, sarkar1} and they have obtained some beautiful results.\\\\
Motivated by above mentioned works, in this paper, we consider conformal $\eta$-Ricci soliton and gradient almost conformal $\eta$-Ricci soliton in the framework of para-Sasakian manifold. We have organized this paper as follows: in first section we look back on some elementary properties of para-Sasakian manifolds; in later section first we prove that if a para-Sasakian manifold satisfies conformal $\eta$-Ricci soliton then the manifold is $\eta$-Einstein and either the soliton vector field $V$ is Killing or it leaves $\phi$ invariant, secondly we prove that if $V$ is pointwise collinear with the characteristic vector field $\xi$ then $V$ is a constant multiple of $\xi$. In the next section, we examine para-Sasakian manifold of dimension greater than 3 with an almost conformal $\eta$ Ricci soliton and show that the manifold is $\eta$-Einstein. Then, we think about a gradient almost conformal $\eta$ Ricci soliton and deduce that the manifold is Einstein and finally, we provide some examples to verify our results.

\section{\textbf{Notes on para-Sasakian manifold}}
A $(2n+1)$-dimensional smooth manifold $M$ is said to have an almost paracontact structure
if it admits a vector field $\xi$, (1,1)-tensor field $\phi$ and a 1-form $\eta$ satisfying the following conditions
\begin{flalign}
    i) \phi^2 = I-\eta\otimes\xi, \label{B1} &&
\end{flalign}
\begin{flalign}
    ii) \eta(\xi) = 1,\label{B2}&&
\end{flalign}
$iii)$ $\phi$ induces on the 2n-dimensional distribution $\mathcal{D}\equiv ker(\eta)$, an almost paracomplex structure $\mathcal{P}$ i.e., $\mathcal{P}^2\equiv I_{\chi(M)}$ and the eigensubbundles $\mathcal{D}^+$ and $\mathcal{D}^-$, corresponding to the eigenvalues $1$, $-1$ of $\mathcal{P}$ respectively, have equal dimension $n$; hence $\mathcal{D}=\mathcal{D}^+\oplus\mathcal{D}^-$.\\
The vector field $\xi$ is called characteristic vector field or Reeb vector field. An immediate consequence of those relations are
\begin{eqnarray}
  \phi\xi &=& 0 \\ \label{B3}
  \eta\circ\phi &=& 0.\label{B4}
\end{eqnarray}
The tensor field $\phi$ induces an almost paracomplex structure on each fibre of $Ker(\eta)$ i.e., the eigendistributions corresponding to eigenvalues $1$ and $-1$ have same dimension $n$. A pseudo-Riemannian metric $g$ is said to be compatible with the almost paracontact structure if
\begin{equation}\label{B5}
  g(\phi X,\phi Y)=-g(X,Y)+\eta(X)\eta(Y)
\end{equation}
holds for arbitrary vector fields $X$ and $Y$ and $(M,\phi,\xi,\eta,g)$ is called an almost paracontact metric manifold.\par
On an almost paracontact metric manifold fundamental 2-form $\Phi$ is defined by $\Phi(X,Y)=g(X,\phi Y)$ for all vector fields $X$ and $Y$ on $M$. An almost paracontact metric manifold for which
\begin{equation}\label{B6}
  \Phi(X,Y)=d\eta(X,Y)=g(X,\phi Y)
\end{equation}
is said to be paracontact metric manifold. In this case, $\eta$ becomes a contact form i.e., $\eta\wedge(d\eta)^n\neq0$ and the manifold becomes a contact manifold. On a paracontact metric manifold $M^{2n+1}(\phi,\xi,\eta,g)$ we consider a self-adjoint operator $h=\frac{1}{2}\mathcal{L}_\xi\phi$, where $\mathcal{L}_\xi$ denotes the Lie derivative along $\xi$. This operator $h$ is symmetric and satisfies
\begin{align*}
  h\phi &= -\phi h, & h\xi &= 0, & \nabla_X\xi &= -\phi X+\phi hX,
\end{align*}
where $\nabla$ is the operator of covariant differentiation w.r.t. the metric $g$. The normality of a paracontact metric manifold $(M,\phi,\xi,\eta,g)$ is equivalent to vanishing of the $(1,2)$-torsion tensor defined by $N_\phi(X,Y)=[\phi,\phi](X,Y)-2d\eta(X,Y)\xi$, where $[\phi,\phi](X,Y)=\phi^2[X,Y]+[\phi X,\phi Y]-\phi[X,\phi Y]-\phi[\phi X,Y]$ for any $X,Y\in\chi(M)$. A normal paracontact metric manifold is called a para-Sasakian metric manifold. It is equivalent to say, an almost paracontact metric manifold is called a para-Sasakian manifold if it satisfies
\begin{equation}\label{B7}
  (\nabla_X\phi)Y=-g(X,Y)\xi+\eta(Y)X
\end{equation}
for arbitrary $X,Y\in\chi(M)$. In a para-Sasakian manifold the operator $h$ vanishes and the manifold satisfies,
\begin{eqnarray}
  \nabla_X\xi &=& -\phi X, \label{B8} \\
  R(X,Y)\xi &=& \eta(X)Y-\eta(Y)X, \label{B9}\\
  R(X,\xi)Y &=& g(X,Y)\xi-\eta(Y)X, \label{B10}\\
  Q\xi &=& -2n\xi, \label{B11}
 \end{eqnarray}
for all vector fields $X$ and $Y$ on $M$ and $R$, $Q$ denote Riemannian curvature tensor and Ricci operator associated with the Ricci tensor $S$ defined by $S(X,Y )=g(QX,Y)$.\par
In \cite{Naik}, authors have prove that in a para-Sasakian manifold the following result holds (for proof see lemma-4)
\begin{equation}\label{B12}
  Q\phi = \phi Q.
\end{equation}\par
Prakasha and Veeresha in \cite{prakasha} established another beautiful result on para-Sasakian manifold (see lemma-1) which using (\ref{B12}) can be restated as
 \begin{eqnarray}
  (\nabla_\xi Q)X &=& 0, \label{B13}\\
  (\nabla_X Q)\xi &=& Q\phi X+2n\phi X.\label{B14}
\end{eqnarray}

\section{\textbf{On Conformal $\eta$-Ricci soliton}}
In this section, we have studied conformal $\eta$-Ricci soliton on Para-Sasakian manifold. First we prove the following lemma which has been used to prove the next theorems.
\begin{lem}
If the metric $g$ of a para-Sasakian manifold represents a conformal $\eta$-Ricci soliton, then
\begin{equation}\label{C1}
  \eta(\mathcal{L}_V\xi)=-(\mathcal{L}_V\eta)\xi=\lambda-\frac{p}{2}-\frac{1}{2n+1}-2n+\mu.
\end{equation}
\begin{proof}
As the metric $g$ satisfies conformal $\eta$-Ricci soliton equation (\ref{A1}), using (\ref{B11}), we can easily obtain
  \begin{equation*}
    (\mathcal{L}_Vg)(X,\xi)+2(\lambda-\frac{p}{2}-\frac{1}{2n+1}-2n+\mu)\eta(X)=0
  \end{equation*}
for arbitrary vector field $X$. Lie differentiation of the relation $\eta(X)=g(X,\xi)$ along the soliton vector field $V$ yields $(\mathcal{L}_Vg)(X,\xi)=(\mathcal{L}_V\eta)X-g(X,\mathcal{L}_V\xi)$. Using this in the foregoing equation, we have
\begin{equation}\label{C2}
(\mathcal{L}_V\eta)X-g(X,\mathcal{L}_V\xi)=-2(\lambda-\frac{p}{2}-\frac{1}{2n+1}-2n+\mu)\eta(X).
\end{equation}
Finally taking Lie derivative of (\ref{B2}) along $V$ into account, we can easily obtain our desired result (\ref{C1}).
\end{proof}
\end{lem}

\begin{thm}
Let $M^{2n+1}(\phi,\xi,\eta,g)$ be a para-Sasakian manifold. If the metric $g$ represents a conformal $\eta$-Ricci soliton then the manifold is $\eta$-Einstein and either the soliton vector field $V$ is Killing or it leaves $\phi$ invariant.
\begin{proof}
  Taking covariant derivative of (\ref{B2}) along arbitrary vector field $Y$ and using (\ref{B8}) we can easily have $(\nabla_Y\eta)X=g(\phi X,Y)$.\par
  Since the metric $g$ of the manifold represents a conformal $\eta$-Ricci soliton, taking covariant derivative of (\ref{A1}) along arbitrary vector field $Z$, we obtain
  \begin{equation}\label{C3}
    (\nabla_Z\mathcal{L}_Vg)(X,Y)=-2(\nabla_ZS)(X,Y)-2\mu[g(\phi X,Z)\eta(Y)+g(\phi Y,Z)\eta(X)]
  \end{equation}
  for all vector fields $X$, $Y$ and $Z$ on $M$. Again from Yano\cite{Yano}, we have the following commutation formula
    \begin{eqnarray*}
      (\mathcal{L}_V\nabla_Xg-\nabla_X\mathcal{L}_Vg-\nabla_{[V,X]}g)(Y,Z)&=&-g((\mathcal{L}_V\nabla)(X,Y),Z)\nonumber\\
      &&-g((\mathcal{L}_V\nabla)(X,Z),Y),
    \end{eqnarray*}
    where $g$ is the metric connection i.e., $\nabla g=0$. So the above equation reduces to
    \begin{equation}\label{C4}
      (\nabla_X\mathcal{L}_Vg)(Y,Z)=g((\mathcal{L}_V\nabla)(X,Y),Z)+g((\mathcal{L}_V\nabla)(X,Z),Y)
    \end{equation}
    for all vector fields $X$, $Y$, $Z$ on $M$. Combining (\ref{C3}) and (\ref{C4}), we have
    \begin{align*}
     g((\mathcal{L}_V\nabla)(X,Z),Y)+&g((\mathcal{L}_V\nabla)(Y,Z),X)=-2(\nabla_ZS)(X,Y)\\
     &-2\mu[g(\phi X,Z)\eta(Y)+g(\phi Y,Z)\eta(X)].
    \end{align*}
     By a straightforward combinatorial computation, the foregoing equation yields
     \begin{eqnarray}\label{C5}
     g((\mathcal{L}_V\nabla)(X,Y),Z) &=& (\nabla_ZS)(X,Y)-(\nabla_XS)(Y,Z)-(\nabla_YS)(Z,X)\nonumber\\
     &&+2\mu[g(\phi X,Z)\eta(Y)+g(\phi Y,Z)\eta(X)]
    \end{eqnarray}
    for all $X,Y,Z\in\chi(M)$. Setting $Y=\xi$ and making use of $(\nabla_ZS)(X,Y)=g((\nabla_ZQ)X,Y)$, (\ref{B12}), (\ref{B13}) and (\ref{B14}), we obtain
    \begin{equation}\label{C6}
      (\mathcal{L}_V\nabla)(X,\xi)=2(\mu-2n)(\phi X)-2Q\phi X.
    \end{equation}
    Differentiating the last equation covariantly with respect to arbitrary vector field $Y$ and using (\ref{B7}) and (\ref{B8}), we acquire
    \begin{eqnarray}\label{C7}
      (\nabla_Y\mathcal{L}_V\nabla)(X,\xi)&=&(\mathcal{L}_V\nabla)(X,\phi Y)-2(\nabla_YQ)(\phi X)-2\eta(X)(QY)\nonumber\\
      &&-2\mu g(X,Y)\xi+2(\mu-2n)\eta(X)Y.
    \end{eqnarray}
    From Yano \cite{Yano}, we know $(\mathcal{L}_VR)(X,Y)\xi=(\nabla_X\mathcal{L}_V\nabla)(Y,\xi)-(\nabla_Y\mathcal{L}_V\nabla)(X,\xi)$. Substituting the values from (\ref{C7}) and using (\ref{B13}), (\ref{C6}), we get
    \begin{equation}\label{C8}
      (\mathcal{L}_VR)(X,\xi)\xi=4(\mu-2n)X-4QX-4\mu\eta(X)\xi,
    \end{equation}
    which holds for an arbitrary vector field $X$. Again, from (\ref{B9}) we get $R(X,\xi)\xi=\eta(X)\xi-X$. By virtue of (\ref{C1}) and (\ref{C2}), Lie differentiation of the last relation along $V$ yields
    \begin{equation}\label{C9}
      (\mathcal{L}_VR)(X,\xi)\xi=2(\lambda-\frac{p}{2}-\frac{1}{2n+1}-2n+\mu)[X-\eta(X)\xi].
    \end{equation}
    Substituting the value of $(\mathcal{L}_VR)(X,\xi)\xi$ from (\ref{C8}) in the foregoing equation and taking inner product with arbitrary vector field $Y$, we obtain
    \begin{eqnarray}\label{C10}
      S(X,Y) &=& \frac{1}{2}(\lambda-\frac{p}{2}-\frac{1}{2n+1}-2n-\mu)\eta(X)\eta(Y) \nonumber\\
      && -\frac{1}{2}(\lambda-\frac{p}{2}-\frac{1}{2n+1}+2n-\mu)g(X,Y)
    \end{eqnarray}
    for all $X,Y\in\chi(M)$. This transforms the soliton equation (\ref{A1}) to
    \begin{equation}\label{C11}
      (\mathcal{L}_Vg)(X,Y)=(\frac{p}{2}+\frac{1}{2n+1}+2n-\lambda-\mu)[g(X,Y)+\eta(X)\eta(Y)].
    \end{equation}
    Differentiating (\ref{C10}) covariantly along arbitrary vector field $Z$, we get \\ $(\nabla_ZS)(X,Y)=\frac{1}{2}(\lambda-\frac{p}{2}-\frac{1}{2n+1}-2n-\mu)[g(\phi X,Z)\eta(Y)+g(\phi Y,Z)\eta(X)]$. Repeated use of this in (\ref{C3}), gives rise to
    \begin{equation}\label{C12}
      (\mathcal{L}_V\nabla)(X,Y)=(\lambda-\frac{p}{2}-\frac{1}{2n+1}-2n+\mu)[\eta(Y)(\phi X)+\eta(X)(\phi Y)]
    \end{equation}
    for arbitrary vector fields $X$ and $Y$ on $M$. Covariant differentiation of the aforementioned equation along arbitrary vector field $Z$ and use of (\ref{B7}), yields
    \begin{align*}
      (\nabla_Z\mathcal{L}_V\nabla)(X,Y) =& (\lambda-\frac{p}{2}-\frac{1}{2n+1}-2n+\mu)[g(\phi Y,Z)(\phi X)+g(\phi X,Z) \\
      & (\phi Y)-g(X,Z)\eta(Y)\xi-g(Y,Z)\eta(X)\xi+2\eta(X)\eta(Y)Z].
    \end{align*}
    Using this relation in $(\mathcal{L}_VR)(X,Y)Z=(\nabla_X\mathcal{L}_V\nabla)(Y,Z)-(\nabla_Y\mathcal{L}_V\nabla)(X,Z)$ (for details see Yano \cite{Yano}) and contracting $X$, we get
    \begin{equation}\label{C13}
      (\mathcal{L}_VS)(Y,Z)=2(\lambda-\frac{p}{2}-\frac{1}{2n+1}-2n+\mu)[(2n+1)\eta(Y)\eta(Z)-g(Y,Z)].
    \end{equation}
    Taking Lie derivative of (\ref{C10}) along $V$ and using (\ref{C11}), we achieve
    \begin{align}\label{C14}
      (\mathcal{L}_VS)(Y,Z) =& \frac{1}{2}(\lambda-\frac{p}{2}-\frac{1}{2n+1}-2n-\mu)[((\mathcal{L}_V\eta)Y)\eta(Z)+ \nonumber \\
      & \eta(Y)((\mathcal{L}_V\eta)Z)]+\frac{1}{2}(\lambda-\frac{p}{2}-\frac{1}{2n+1}+2n-\mu) \nonumber\\
      &(\lambda-\frac{p}{2}-\frac{1}{2n+1}-2n+\mu)[g(Y,Z)+\eta(Y)\eta(Z)].
    \end{align}
    Comparisons of  (\ref{C13}) and (\ref{C14}) gives
    \begin{align}\label{C15}
       &2(\lambda-\frac{p}{2}-\frac{1}{2n+1}-2n+\mu)[(2n+1)\eta(Y)\eta(Z)-g(Y,Z)]=\nonumber\\
       & \frac{1}{2}(\lambda-\frac{p}{2}-\frac{1}{2n+1}-2n-\mu)[((\mathcal{L}_V\eta)Y)\eta(Z)+\eta(Y)((\mathcal{L}_V\eta)Z)] \nonumber \\
      & +\frac{1}{2}(\lambda-\frac{p}{2}-\frac{1}{2n+1}+2n-\mu)(\lambda-\frac{p}{2}-\frac{1}{2n+1}-2n+\mu) \nonumber\\
      &[g(Y,Z)+\eta(Y)\eta(Z)].
    \end{align}
    Substituting $Y$ and $Z$ by $\phi^2Y$ and $\phi Z$ respectively and using (\ref{B1}), (\ref{B4}) and (\ref{B6}), we obtain
    \begin{equation}\label{C16}
      (\lambda-\frac{p}{2}-\frac{1}{2n+1}-2n+\mu)(\lambda-\frac{p}{2}-\frac{1}{2n+1}+2n-\mu+4)d\eta(Y,Z)=0
    \end{equation}
    $\forall Y,Z\in\chi(M)$. As we know, in para-Sasakian manifold $d\eta\neq0$, we have $(\lambda-\frac{p}{2}-\frac{1}{2n+1}-2n+\mu)(\lambda-\frac{p}{2}-\frac{1}{2n+1}+2n-\mu+4)=0$. This gives either $\lambda=\frac{p}{2}+\frac{1}{2n+1}+2n-\mu$ or $\lambda=\frac{p}{2}+\frac{1}{2n+1}-2n+\mu-4$.\par
    \medskip
    \textbf{Case-I:} If $\lambda=\frac{p}{2}+\frac{1}{2n+1}+2n-\mu$, then (\ref{C10}) reduces to $S(X,Y)=(\mu-2n)g(X,Y)-\mu\eta(X)\eta(Y)$ i.e., the manifold is $\eta$-Einstein. Also (\ref{C11}) gives $\mathcal{L}_Vg=0$. So, $V$ is Killing vector field.\par
    \medskip
    \textbf{Case-II:} Using $\lambda=\frac{p}{2}+\frac{1}{2n+1}-2n+\mu-4$ in (\ref{C10}), we get $S(X,Y)=2g(X,Y)-2(n+1)\eta(X)\eta(Y)$. So, the manifold is $\eta$-Einstein. Substituting $Y$ by $\phi Y$ and setting $Z=\xi$ in (\ref{C15}), we obtain $(\mathcal{L}_V\eta)(\phi Y)=0$. Further, replacing $Y$ by $\phi Y$ and using (\ref{B1}), (\ref{C1}) and $\lambda=\frac{p}{2}+\frac{1}{2n+1}-2n+\mu-4$, we get
    \begin{equation*}
      \mathcal{L}_V\eta=2(2n-\mu+2)\eta
    \end{equation*}
    Exterior differentiation of the foregoing equation and use of well-known relation $d(\mathcal{L}_V\eta)=\mathcal{L}_Vd\eta$ and (\ref{B6}), yields
    \begin{equation}\label{C17}
      (\mathcal{L}_Vd\eta)(X,Y)=2(2n-\mu+2)g(X,\phi Y)
    \end{equation}
    for arbitrary vector fields $X$ and $Y$ on $M$. Lie differentiation of (\ref{B6}) along $V$, infers
    \begin{equation}\label{C18}
      (\mathcal{L}_Vd\eta)(X,Y)=2(2n-\mu+2)g(X,\phi Y)+g(X,(\mathcal{L}_V\phi)Y)
    \end{equation}
    $\forall X,Y\in\chi(M)$. Comparing this with (\ref{C17}) gives $\mathcal{L}_V\phi=0$, as $X$ and $Y$ are arbitrary vector fields. So, $V$ leaves $\phi$ invariant.
\end{proof}
\end{thm}

\begin{thm}
If the metric $g$ of a para-Sasakian manifold $M$ represents a conformal $\eta$-Ricci soliton and if the soliton vector field $V$ is pointwise collinear with the characteristic vector field $\xi$ then $V$ is a constant multiple of $\xi$ and the scalar curvature of the manifold is constant.
\begin{proof}
  Since the soliton vector field $V$ is pointwise collinear with the characteristic vector field $\xi$, so, $V=f\xi$ where, $f$ is a smooth function on $M^{2n+1}$. Substituting $V=f\xi$ in $(\mathcal{L}_Vg)(X,Y)=g(\nabla_XV,Y)+g(X,\nabla_YV)$ and using (\ref{B8}), we get
  \begin{equation}\label{D1}
    (\mathcal{L}_Vg)(X,Y)=(Xf)\eta(Y)+(Yf)\eta(X)
  \end{equation}
  for arbitrary vector field $X$ and $Y$ on $M$. Using (\ref{D1}) in the soliton equation (\ref{A1}), we obtain
  \begin{eqnarray}\label{D2}
    &&(Xf)\eta(Y)+(Yf)\eta(X)+2S(X,Y)+2\mu\eta(X)\eta(Y)\nonumber\\
    &&+2(\lambda-\frac{p}{2}-\frac{1}{2n+1})g(X,Y)=0.
  \end{eqnarray}
  Setting $Y=\xi$ and using (\ref{B2}) and (\ref{B11}), the above equation becomes
  \begin{equation}\label{D3}
    Df=2[2n+\frac{p}{2}+\frac{1}{2n+1}-\frac{1}{2}(\xi f)-\lambda-\mu]\xi.
  \end{equation}
  Taking inner product with respect to Reeb vector field $\xi$, we acquire
  \begin{equation}\label{D4}
    \xi f=2n+\frac{p}{2}+\frac{1}{2n+1}-\lambda-\mu.
  \end{equation}
  From previous theorem we obtained either $\lambda=\frac{p}{2}+\frac{1}{2n+1}+2n-\mu$ or $\lambda=\frac{p}{2}+\frac{1}{2n+1}-2n+\mu-4$.
  If we consider $\lambda=\frac{p}{2}+\frac{1}{2n+1}+2n-\mu$, then from (\ref{D4}) we get $\xi f=0$. Substituting these values in (\ref{D3}), yields
  \begin{equation}\label{D5}
    Df=0.
  \end{equation}
  Now, if we consider $\lambda=\frac{p}{2}+\frac{1}{2n+1}-2n+\mu-4$, then from (\ref{D4}) we obtain $\xi f=2(2n-\mu+2)$ and from (\ref{D3}), we get
  \begin{equation}\label{D6}
    Df=(\xi f)\xi.
  \end{equation}
  Taking (\ref{B8}) into account, differentiating the foregoing equation along $X$ and then taking scalar product with arbitrary vector field $Y$, leads to
  \begin{equation}\label{D7}
    g(\nabla_XDf,Y)=(X(\xi f))\eta(Y)-(\xi f)g(\phi X,Y).
  \end{equation}
  Anti-symmetrizing the last equation and using $g(\nabla_XDf,Y)=g(X,\nabla_YDf)$, we have
  \begin{equation}\label{D8}
    (X(\xi f))\eta(Y)-(Y(\xi f))\eta(X)-2(\xi f)g(\phi X,Y)=0
  \end{equation}
  for arbitrary vector fields $X$ and $Y$. If we let $X$ to be unit vector (i.e., $g(X,X)=1$) in $Ker(\eta)$, then $\phi X$ also becomes a unit vector with $g(\phi X,\phi Y)=-1$. Now replacing $Y$ by $\phi X$ in (\ref{D8}), we get $\xi f=0$. Substituting this value in (\ref{D6}) leads to
  \begin{equation}\label{D9}
    Df=0.
  \end{equation}
  Combining (\ref{D5}) and (\ref{D9}) we can conclude that $Df=0$ in entire manifold. Therefore $f$ is constant. So, $V$ is a constant multiple of $\xi$.\\
  The equation (\ref{D2}) reduces to $QX+2(\lambda-\frac{p}{2}-\frac{1}{2n+1})X+2\mu\eta(X)\xi=0$. Tracing of this equation leads to $r=2-2\mu+(2n+1)(p-2\lambda)$, where $r$ denotes the scalar curvature of the manifold. This completes the proof.
\end{proof}
\end{thm}

\section{\textbf{On almost conformal $\eta$-Ricci soliton}}
In this section, we consider almost conformal $\eta$-Ricci soliton on para-Sasakian manifold. It follows from (\ref{A1}) that almost conformal $\eta$-Ricci soliton is the generalization of almost $\eta$-Ricci soliton because it involve two smooth functions $\lambda$ and $\mu$.
\begin{defn}
A vector field $V$ is said to be an infinitesimal contact transformation if there exists a certain $a\in C^\infty(M)$ such that
\begin{equation}\label{E1}
  \mathcal{L}_V\eta=a\eta.
\end{equation}
If $a=0$, then $V$ is called strictly infinitesimal contact transformation.
\end{defn}

\begin{thm}
Let $M^{2n+1}$ be a para-Sasakian manifold with $n>1$. If $g$ represents an almost conformal $\eta$-Ricci soliton with the soliton vector field $V$ as infinitesimal contact transformation, then the manifold is $\eta$-Einstein and either the soliton vector field $V$ is Killing or it leaves $\phi$ invariant.
\begin{proof}
  Taking $d(\mathcal{L}_V\eta)=\mathcal{L}_Vd\eta$ into account, exterior derivative of (\ref{E1}) gives
  \begin{equation*}
    \mathcal{L}_Vd\eta=(da)\wedge\eta+ad\eta.
  \end{equation*}
  Using (\ref{B6}), the foregoing equation can be rewritten as
  \begin{equation}\label{E2}
    (\mathcal{L}_Vd\eta)(X,Y)=\frac{1}{2}[(Xa)\eta(Y)-\eta(X)(Ya)]+ag(X,\phi Y)
  \end{equation}
  for arbitrary vector fields $X$ and $Y$ on $M$. Lie differentiation of (\ref{B6}) along the soliton vector field $V$, yields
  \begin{equation}\label{E3}
    (\mathcal{L}_Vd\eta)(X,Y)=g(X,(\mathcal{L}_V\phi)Y)-2g((\lambda-\frac{p}{2}-\frac{1}{2n+1})X+QX,\phi Y)
  \end{equation}
  $\forall X,Y\in\chi(M)$. Comparing the aforementioned equation with (\ref{E2}), we obtain
  \begin{equation}\label{E4}
    2(\mathcal{L}_V\phi)Y=\eta(Y)(Da)-(Ya)\xi+(2a+4\lambda-2p-\frac{4}{2n+1})(\phi Y)+4Q\phi Y
  \end{equation}
  for any vector $Y$ in $\chi(M)$. Setting $Y=\xi$ and using (\ref{B2}) and (\ref{B3}), we get
  \begin{equation}\label{E5}
    2(\mathcal{L}_V\phi)\xi=Da-(\xi a)\xi.
  \end{equation}
  Taking (\ref{A1}) and (\ref{E1}) into consideration, Lie differentiating  $g(\xi,\xi)=1$ and $g(X,\xi)=\eta(X)$ along $V$, we achieve
  \begin{eqnarray}
    \eta(\mathcal{L}_V\xi) &=& \lambda-\frac{p}{2}-\frac{1}{2n+1}-2n+\mu \label{E6}\\
    \mathcal{L}_V\xi &=& (a+2\lambda-p-\frac{2}{2n+1}-4n+2\mu)\xi. \label{E7}
  \end{eqnarray}
  Combining above two relations, we get
  \begin{equation}\label{E8}
    a=2n-\lambda-\mu+\frac{p}{2}+\frac{1}{2n+1}.
  \end{equation}
  Lie differentiating (\ref{B3}) along the soliton vector field $V$ and using (\ref{E8}), we obtain $(\mathcal{L}_V\phi)\xi=0$. Combining this with (\ref{E5}) we have $Da=(\xi a)\xi$, which further implies
  \begin{equation}\label{E9}
    da=(\xi a)\eta.
  \end{equation}
  Operating the foregoing equation by exterior derivative operator $d$ and using $d^2=0$, yields
  \begin{equation*}
    d(\xi a)\wedge\eta+(\xi a)d\eta=0.
  \end{equation*}
  Taking $\eta\wedge\eta=0$ and $\eta\wedge d\eta\neq 0$ into account, wedge product of the above equation with respect to the 1-form $\eta$ gives $\xi a=0$. Substituting this in (\ref{E9}) we get $da=0$ and so $a$ is constant. Then the equation (\ref{E4}) reduces to
  \begin{equation}\label{E10}
    (\mathcal{L}_V\phi)Y=(a+2\lambda-p-\frac{2}{2n+1})(\phi Y)+2Q\phi Y
  \end{equation}
  for any vector field $Y$ on $M$. Operating (\ref{B1}) by $\mathcal{L}_V$ and using (\ref{E1}) and (\ref{E7}) we get $\mathcal{L}_V\phi^2=0$. Consequently, we get $(\mathcal{L}_V\phi)(\phi X)+\phi(\mathcal{L}_V\phi)X=0$ for an arbitrary vector field $X$. After repeated application of (\ref{E10}) and use of (\ref{B1}), (\ref{B11}) and (\ref{B12}), the aforementioned relation leads to
  \begin{equation}\label{E11}
    QX=-(\frac{a}{2}+\lambda-\frac{p}{2}-\frac{1}{2n+1})X+(\frac{a}{2}+\lambda-\frac{p}{2}-\frac{1}{2n+1}-2n)\eta(X)\xi
  \end{equation}
  for all field $X$ of $M$. Covariant derivative of (\ref{E11}) along an aritrary vector field $Y$ yields
  \begin{align}\label{E12}
    (\nabla_YQ)X=&(\frac{a}{2}+\lambda-\frac{p}{2}-\frac{1}{2n+1}-2n)[g(Y,\phi X)\xi-\eta(X)(\phi Y)] \nonumber\\
    & -(Y\lambda)[X-\eta(X)\xi].
  \end{align}
  It is well-known that $Xr=g((\nabla_XQ)e_i,e_i)$ and $\frac{1}{2} Xr=g((\nabla_{e_i}Q)X,e_i)$, where $\{e_i\}_{i=1}^{2n+1}$ is an orthonormal basis of the manifold. Successive application of (\ref{E12}) in these two relations infers
  \begin{eqnarray}
    Xr &=& -2n(X\lambda), \label{E13}\\
    Xr &=& -2(X\lambda)+2\eta(X)(\xi\lambda) \label{E14}
  \end{eqnarray}
  for an arbitrary vector field $X$. Since the characteristic vector field $\xi$ is a Killing vector field in a para-Sasakian manifold, it follows that $\xi r=0$. Plugging this in (\ref{E13}), we obtain $\xi\lambda=0$ as we assume $n\neq 0$. Consequently, (\ref{E14}) reduces to $Xr=-2(X\lambda), \forall X\in\chi(M)$. Setting this in (\ref{E13}), we get
  \begin{equation*}
    (n-1)(X\lambda)=0.
  \end{equation*}
  Since $n>1$ and $X$ is an arbitrary vector field, we conclude that $\lambda$ is a constant. Thus it follows from (\ref{E8}) that $\mu$ is also constant. Then the soliton reduces to conformal $\eta$-Ricci soliton and the result follows from theorem-3.2.
\end{proof}
\end{thm}

\begin{thm}
Let $M^{2n+1}$ be a para-Sasakian manifold of dimension $>3$. If $g$ represents a gradient almost conformal $\eta$-Ricci soliton, then the soliton reduces to gradient almost conformal Ricci soliton and the manifold is Einstein.
\begin{proof}
  From the gradient almost conformal $\eta$-Ricci soliton equation (\ref{A2}), we easily obtain
  \begin{equation}\label{F1}
    \nabla_XDf=-QX+(\frac{p}{2}+\frac{1}{2n+1}-\lambda)X-\mu\eta(X)\xi
  \end{equation}
  for an arbitrary vector field $X$ of $M$. Taking covariant derivative along an arbitrary vector field $Y$, we acquire
  \begin{align*}
    \nabla_Y\nabla_XDf=& -(\nabla_YQ)X-Q(\nabla_YX)+(\frac{p}{2}+\frac{1}{2n+1}-\lambda)(\nabla_YX)-(Y\lambda)X\\
    & -(Y\mu)\eta(X)\xi-\mu[g(\phi X,Y)\xi+\eta(\nabla_YX)\xi-\eta(X)(\phi Y)].
  \end{align*}
  From Yano \cite{Yano}, we know $R(X,Y)Z=\nabla_X\nabla_YZ-\nabla_Y\nabla_XZ-\nabla_{[X,Y]}Z$ for all $X,Y,Z\in\chi(M)$. Plugging the above equation along with (\ref{F1}) in this curvature property, we infer
  \begin{align}\label{F2}
    R(X,Y)Df &= (\nabla_YQ)X-(\nabla_XQ)Y+(Y\lambda)X-(X\lambda)Y+(Y\mu)\eta(X)\xi \nonumber\\
    & -(X\mu)\eta(Y)\xi+\mu[2g(\phi X,Y)\xi-\eta(X)(\phi Y)+\eta(Y)(\phi X)]
  \end{align}
  for all $X,Y\in\chi(M)$. Setting $Y=\xi$ in the foregoing equation and using (\ref{B13}) and (\ref{B14}), we get
  \begin{equation*}
    R(X,\xi)Df=-Q\phi X+(\mu-2n)(\phi X)+(\xi\lambda)X-(X\lambda)\xi+(\xi\mu)\eta(X)\xi-(X\mu)\xi.
  \end{equation*}
  Combining the last equation with (\ref{B7}) and (\ref{B10}), we obtain
  \begin{align}\label{F3}
    g((\nabla_X\phi)Y,Df)= & -g(Q\phi X,Y)+(\mu-2n)g(\phi X,Y)+(\xi\lambda)g(X,Y) \nonumber\\
    & -(X\lambda)\eta(Y)+(\xi\mu)\eta(X)\eta(Y)-(X\mu)\eta(Y)
  \end{align}
  for arbitrary vector fields $X$ and $Y$ on $M$. Replacing $X$ and $Y$ by $\phi X$ and $\phi Y$ in the foregoing equation, we achieve
  \begin{align}\label{F4}
    g((\nabla_{\phi X}\phi)\phi Y,Df)= & g(Q\phi X,Y)-(\mu-2n)g(\phi X,Y) \nonumber\\
    &-(\xi\lambda)[g(X,Y)-\eta(X)\eta(Y)]
  \end{align}
  where we have used (\ref{B1}), (\ref{B2}), (\ref{B5}), (\ref{B11}) and (\ref{B12}). Subtraction of (\ref{F3}) from (\ref{F4}) yields
  \begin{align}\label{F5}
    &g((\nabla_{\phi X}\phi)\phi Y-(\nabla_X\phi)Y,Df)= 2g(Q\phi X,Y)-2(\mu-2n)g(\phi X,Y)- \nonumber\\
    &2(\xi\lambda)g(X,Y)+(\xi\lambda)\eta(X)\eta(Y)+(X\lambda)\eta(Y)-(\xi\mu)\eta(X)\eta(Y)+(X\mu)\eta(Y).
  \end{align}
  From Zamkovoy \cite{Zamkovoy}, we know
  \begin{equation}\label{F6}
    (\nabla_{\phi X}\phi)\phi Y-(\nabla_X\phi)Y=2g(X,Y)\xi-\eta(Y)[X+\eta(X)\xi]
  \end{equation}
  holds for arbitrary vector fields $X$ and $Y$ in a para-Sasakian manifold (for proof see lemma-2.7 of \cite{Zamkovoy} and here we have used $h=0$ which holds in para-Sasakian manifold). Taking (\ref{F6}) into account, (\ref{F5}) can be rewritten as
  \begin{align}\label{F7}
    & 2g(X,Y)(\xi f)-\eta(Y)(Xf)-\eta(X)\eta(Y)(\xi f)= 2g(Q\phi X,Y) \nonumber\\
    & -2(\mu-2n)g(\phi X,Y)-2(\xi\lambda)g(X,Y)+(\xi\lambda)\eta(X)\eta(Y)+\nonumber\\
    & (X\lambda)\eta(Y)-(\xi\mu)\eta(X)\eta(Y)+(X\mu)\eta(Y).
  \end{align}
  Anti-symmetrizing the last equation and then replacing $X$ and $Y$ by $\phi X$ and $\phi Y$, respectively we get $Q\phi X=(\mu-2n)(\phi X)$. Further, substitution of $X$ by $\phi X$ in the last relation yields
  \begin{equation}\label{F8}
    QX=(\mu-2n)X-\mu\eta(X)\xi.
  \end{equation}
  Covariant differentiation of (\ref{F8}) along an arbitrary vector field $Y$ and using that expression of $(\nabla_YQ)X$ in (\ref{F2}), we obtain
  \begin{equation}\label{F9}
    R(X,Y)Df= (Y\mu)X-(X\mu)Y+(Y\lambda)X-(X\lambda)Y
  \end{equation}
  $\forall X,Y\in\chi(M)$. Contraction of (\ref{F2}) and (\ref{F9}) over $X$ yields
  \begin{eqnarray}
    Q(Df) &=& \frac{1}{2}(Dr)+2n(D\lambda)+(D\mu)-(\xi\mu)\xi \label{F10}\\
    Q(Df) &=& 2n(D\mu+D\lambda) \label{F11}
  \end{eqnarray}
  Comparing the last two relations we get
  \begin{equation}\label{F12}
    (2n-1)D\mu=\frac{1}{2}(Dr)-(\xi\mu)\xi.
  \end{equation}
  Tracing (\ref{F8}), we find $r=2n\mu-2n(2n+1)$. So, $Dr=2nD\mu$. Plugging this relation in (\ref{F12}), we obtain
  \begin{equation}\label{F13}
    (n-1)D\mu+(\xi\mu)\xi=0.
  \end{equation}
  As we know $g((\nabla_{e_i}Q)\xi,e_i)=\frac{1}{2}(\xi r)$, we can easily obtain $\xi r=0$. Combining this result with $Dr=2nD\mu$, we get $\xi\mu=0$. Using this in (\ref{F13}), we have $D\mu=0$ (since $n\neq 1$), so $\mu$ is constant. Scalar product of (\ref{F11}) with respect to the characteristic vector field $\xi$ and use of (\ref{B11}), yields
  \begin{equation}\label{F14}
    \xi(f+\lambda)=0
  \end{equation}
  as $n>1$. Now, setting $Y=\xi$ in (\ref{F9}) and using (\ref{B10}), we obtain
  \begin{equation}\label{F15}
    (X(f+\lambda))\xi=(\xi(f+\lambda))X
  \end{equation}
  for an arbitrary vector field $X$. Plugging equation (\ref{F14}) in the above equation, we get that $f+\lambda$ is constant. As we have $\mu$ and $f+\lambda$ are constant, replacing $X$ by $Df$ in (\ref{F8}), we get
  \begin{equation}\label{F16}
    \mu[(Df)-(\xi f)\xi]=0.
  \end{equation}
  We suppose, $\mu$ is a non-zero constant. Then, from last equation, we get $Df=(\xi f)\xi$. Substitution of $Df$ by $(\xi f)\xi$ in (\ref{F1}) and use of (\ref{B8}), yields
  \begin{equation}\label{F17}
    (X(\xi f))\xi-(\xi f)(\phi X)=(\frac{p}{2}+\frac{1}{2n+1}-\lambda-\mu+2n)X
  \end{equation}
  for any vector field $X$ of $\chi(M)$. Taking inner product of the foregoing equation with respect to $\xi$ and substituting the resultant relation in (\ref{F17}), we get
  \begin{equation*}
    (\xi f)(\phi X)+(\frac{p}{2}+\frac{1}{2n+1}-\lambda-\mu+2n)[X-\eta(X)\xi]=0.
  \end{equation*}
  Contracting the last equation over $X$ and using $n\neq 0$, we obtain $\frac{p}{2}+\frac{1}{2n+1}-\lambda-\mu+2n=0$. Since $\mu$ is a constant, we have $\lambda$ is also constant. As we know $f+\lambda$ is constant, this yields $f$ is constant. This contradicts the fact that the soliton vector field $V$ is non-zero as we get $V=Df=0$.\\
  So, $\mu$ must be identically equal to zero and the soliton reduces to gradient almost conformal Ricci soliton. Finally, the equation (\ref{F8}) reduces to $QX=-2nX~\forall X\in\chi(M)$. So, the manifold becomes Einstein with Einstein constant $-2n$. Furthermore, the scalar curvature of the manifold can be expressed as $r=-2n(2n+1)$.
\end{proof}
\end{thm}
\begin{ex}
We consider the example of the paper \cite{Naik}. In this paper, authors considered the Euclidean space $M=\Bbb{R}^3$  with Cartesian coordinates $(x, y, z)$ and defined the normal almost paracontact metric structure $(\varphi, \xi, \eta, g)$ on $M$ as follows:
  \begin{align*}
\varphi(\frac{\partial}{\partial x})&=\frac{\partial}{\partial y}, & \varphi(\frac{\partial}{\partial y})&= \frac{\partial}{\partial x}-y\frac{\partial}{\partial z}, & \varphi(\frac{\partial}{\partial z})&=0,
\end{align*}
\begin{align*}
\xi&=2\frac{\partial}{\partial z}, & \eta&=\frac{1}{2}(dz+ydx),
\end{align*}
\begin{equation*}
    (g_{ij})=\begin{pmatrix}
               \frac{y^2 - 1}{2} & 0 & \frac{y}{4}  \\
               0 & \frac{1}{4} & 0  \\
               \frac{y}{4} & 0 & \frac{1}{4}
             \end{pmatrix}
  \end{equation*}
and authors has shown that the manifold $M$ is para-Sasakian. Authors has taken pseudo-orthonormal $\varphi$-basis $e_1=2\partial y$, $e_2=2\partial x-2y\partial z$ and $e_3=\xi=2\partial z$ and also obtained the expressions of the curvature tensor and the
Ricci tensor respectively as follows:
\begin{align*}
  R(e_1, e_2)e_1 &=-3e_2, & R(e_1, e_2)e_2 &=-3e_1, & R(e_1, e_2)e_3 &=0,\\
  R(e_1, e_3)e_1 &=\xi, & R(e_1, e_3)e_2 &=0, & R(e_1, e_3)e_3 &= -e_1,\\
  R(e_2, e_3)e_1 &= 0, & R(e_2, e_3)e_2 &= -\xi, & R(e_2, e_3)e_3 &= -e_2,
\end{align*}
and
\begin{align*}
S(e_1, e_1) &=2, & S(e_2, e_2) &=-2, & S(e_3, e_3)&=-2.
\end{align*}
Also the scalar curvature $r$=2. Thus
\begin{equation}\label{e1}
S(X, Y)=2g(X, Y)-4\eta(X)\eta(Y)~~\forall X,Y\in\chi(M).
\end{equation}
Let $f:M\rightarrow \mathbb{R}$ be a smooth function defined by,
\begin{equation}\label{e3}
  f(x,y,z)=\frac{x^2}{2}+\frac{y^2}{2}+z^2.
\end{equation}
Then the gradient of $f$, $Df$ is given by,
\begin{equation}
Df=(x\frac{\partial}{\partial x}+y\frac{\partial}{\partial y}+2z\frac{\partial}{\partial z})
\end{equation}
Now $(\mathcal{L}_{Df}g)(e_{1}, e_{1})=-2g(e_{1}, e_{1})=2$,  $(\mathcal{L}_{Df}g)(e_{2}, e_{2})=2g(e_{2},e_{2})=-2$ and $(\mathcal{L}_{Df}g)(e_{3}, e_{3})=4$.
Then from the above results we can verify that,
 \begin{equation}\label{e2}
(\mathcal{L}_{Df}g)(X, Y)=2\{g(X, Y)+\eta(X)\eta(Y)\}
\end{equation}
for all $X, Y\in\chi(M).$ From (\ref{e1}) and (\ref{e2}) we obtain that $g$ represents a gradient almost conformal $\eta$-Ricci soliton i.e., it satisfies (\ref{A2}) for $V=Df$, where $f$ is defined by (\ref{e3}), $\lambda=\frac{p}{2}-\frac{8}{3}$ and $\mu=3$.

\end{ex}

\section{\textbf{Conclusion and Remarks}}
The study of conformal Ricci solitons and conformal $\eta$ Ricci solitons on Riemannian manifolds and pseudo-Riemannian manifolds is a great importance in the area of differential geometry, specially in Riemannian geometry and in special relativistic physics as well. Basically, conformal Ricci flow is the most prominent flagship of modern physics.  The conformal $\eta$-Ricci soliton is a new notion not only in the area of differentiable manifold but in the area of mathematical physics, general relativity and quantum cosmology, quantum gravity, Black hole as well.  It deals a geometric and physical applications with relativistic viscous fluid spacetime admitting heat flux and stress, dark and dust fluid general relativistic spacetime, radiation era in general relativistic spacetime. The application of this solitons and our results will not only play an important and significant role in paracontact geometry, but also it has a motivational contribution in mathematical fluid dynamics, thermodynamics etc. We can also find the application of such solitons on some Einstein space time like hyper-generalized quasi-Einstein spacetime, mixed super quasi-Einstein spacetime with the connection of general relativity. There are some questions arise from our article to study further research.\\
(i) Is theorem 3.3 true if we assume soliton vector field $V$ is pointwise collinear with the characteristic vector field $\xi$
?\\
(ii) Is theorem 4.2 true without assuming soliton vector field $V$ as an infnitesimal contact transformation?\\
(iii) Are theorem 4.2 and theorem 4.3 true if we consider the dimension of the manifold is 3?

\section{\textbf{Acknowledgements}}
The first author is the corresponding author and this work was financially supported by UGC Senior Research Fellowship of India, Sr. No. 2061540940. Ref. No:21/06/2015(i)EU-V.

\end{document}